%% file: main.tex
\documentclass[journal,twoside,web]{ieeecolor}
\usepackage{tikz}
\usepackage{tikz-network}
\usepackage{pstricks}
\usetikzlibrary{shapes}
\usepackage{pgfplots}
\usetikzlibrary{matrix}
\usepackage{./custom/generic}
\usepackage{cite}
\usepackage{algorithmic}
\usepackage{graphicx}
\usepackage{textcomp}
\usepackage{bm, amsfonts, amsmath, shortcuts_OPT}
\usepackage{graphicx}
\usepackage{url}

\usepackage[colorlinks=true,
linkcolor=green!25!black,
urlcolor=blue,
citecolor=blue]{hyperref}  

\usepgfplotslibrary{groupplots,dateplot}

\usepackage{epsfig} % for postscript graphics files
\usepackage{cite}
\usepackage{lcsys}
\usepackage{amsthm}
\usepackage{algorithm} 

\let\labelindent\relax
\usepackage{enumitem}

\newtheorem{theorem}{Theorem} 
\newtheorem{lemma}{Lemma}
\newtheorem{proposition}{Proposition}
\newtheorem{corollary}{Corollary}

\newtheorem{assumption}{Assumption}

\begin{document}

\def\BibTeX{{\rm B\kern-.05em{\sc i\kern-.025em b}\kern-.08em
    T\kern-.1667em\lower.7ex\hbox{E}\kern-.125emX}}
\markboth{\journalname, VOL. XX, NO. XX, XXXX 2017}
{Author \MakeLowercase{\textit{et al.}}: Preparation of Papers for IEEE Control Systems Letters (August 2022)}

\title{Optimal Pricing for Linear-Quadratic Games with Nonlinear Interaction Between Agents}

\author{
Jiamin Cai, Chenyue Zhang, \IEEEmembership{Student Member, IEEE}, Hoi-To Wai, \IEEEmembership{Member, IEEE}
\thanks{$^{a}$Department of SEEM, The Chinese University of Hong Kong, Shatin, Hong Kong SAR of China; \texttt{jmcai@link.cuhk.edu.hk,\{czhang,htwai\}@se.cuhk.edu.hk}}
}

\maketitle
\thispagestyle{empty}

%%%%%%%%%%%%%%%%%%%%%%%%%%%%%%%%%%%%%%%%%%%%%%%%%%%%%%%%%%%%%%%%%%%%%%%%%%%%%%%%
\begin{abstract}
This paper studies a class of network games with linear-quadratic payoffs and externalities exerted through a strictly concave interaction function. This class of game is motivated by the diminishing marginal effects with peer influences. We analyze the optimal pricing strategy for this class of network game. First, we prove the existence of a unique Nash Equilibrium (NE). Second, we study the optimal pricing strategy of a monopolist selling a divisible good to agents. We show that the optimal pricing strategy, found by solving a bilevel optimization problem, is strictly better when the monopolist knows the network structure as opposed to the best strategy agnostic to network structure. Numerical experiments demonstrate that in most cases, the maximum revenue is achieved with an asymmetric network. These results contrast with the previously studied case of linear interaction function, where a network-independent price is proven optimal with symmetric networks. Lastly, we describe an efficient algorithm for finding the optimal pricing strategy.  
\end{abstract}


\section{Introduction}
\label{sec:introduction}
\IEEEPARstart{A} key feature of social networks is that an agent's action can influence and be influenced by neighbors' actions. To analyze the behavior of social networks, many works \cite{jackson2015games, zhou2016game, parise2017sensitivity} have considered studying network games as a general framework to capture the interaction among agents and analyze the latter's  actions at \emph{equilibrium}.  
Among others, the linear-quadratic game setting is popular \cite{ballester2006s, candogan2012optimal, bramoulle2014strategic} as it captures essential features of social networks while providing useful analysis.

Pricing or intervention strategy describes the decisions, made by an independent entity, to modulate the agents' actions \cite{demange2017optimal, parise2021analysis, galeotti2009influencing, galeotti2020targeting, maheshwari2022inducing} while achieving a certain goal depending on the agents' equilibrium actions. As an example, in social networks, monopolist sets prices for agents who adjust their consumption of a good while the monopolist maximizes his/her revenue from selling the good. However, the optimal pricing problem is a bilevel optimization problem which can be difficult to analyze. Only few results are found in the literature, e.g., \cite{candogan2012optimal, ballester2006s, galeotti2009influencing} analyzed the optimal pricing strategy for linear-quadratic games and drew connections between the optimal prices and network's centrality, \cite{demange2017optimal} studied the effects of general strategic complementary, \cite{ata2023latent} studied the inverse problem of social network from pricing experiments. 

This paper is concerned with the optimal pricing problem where a monopolist seeks to maximize his/her total revenue made from selling a divisible good to the agents, where the agents are incentivized to buy the good through pricing in a linear quadratic game. We focus on a strategic complementary setting \cite{belhaj2014network} such that an agent is positively influenced by peers through a \emph{strictly concave interaction function}. This is unlike the model of {linear interaction} studied in \cite{candogan2012optimal, ballester2006s} as we are motivated by the \emph{diminishing marginal} effect {of peers \cite{cooley2007desegregation}, e.g., in investment problems as modeled by \cite{zhou2016game}.}

Prior works have studied network games with various forms of nonlinear interaction among agents. Many of these works have shown the uniqueness of Nash Equilibrium (NE) for the game, e.g., \cite{belhaj2014network,naghizadeh2017on} considered concave function on the aggregated influences and extended the criteria for the unique NE to exist as well as studying the property of NE, \cite{shende2021network} focused on linear best response with nonlinear social warfare, \cite{wang2023network, ebrahimi2023united} considered the equilibrium in multiplex networks. These works are related to studies on network aggregative games \cite{scutari2014real, parise2019variational}, where distributed algorithms for computing the NE are in \cite{salehisadaghiani2016distributed, liang2017distributed, tatarenko2020geometric}.
For intervention strategy, the closest work to ours is \cite{demange2017optimal}, which shows that the targeting strategy optimal to an objective of aggregative action is affected by concavity/convexity of reaction function.
{Our work is also related to studies on Stackelberg equilibrium \cite{bacsar1998dynamic}. Recent works have focused on the computational aspects \cite{briest2012stackelberg,bohnlein2021revenue}.}

Motivated by the above works, we study the effects on the optimal pricing problem caused by the nonlinear interaction function that describes the diminishing marginal phenomena.
We provide a few new analytical results: 
\begin{itemize}[leftmargin=*]
    \item When given a pricing strategy, we first prove a sufficient condition for the existence of unique NE. Our condition is characterized by the Lipschitzness of interaction function.
    \item For optimal pricing problem, we prove that the optimal total revenue is \emph{strictly better} than the revenue from a \emph{network agnostic} strategy found without network information. Importantly, we bound the \emph{price of information (PoI)}, defined as the ratio between the total revenues made by the optimal and network agnostic strategy. We show that the \emph{PoI} increases with the curvature of interaction function.
    \item We show that the optimal pricing problem is equivalent to a strongly concave optimization problem through reparameterizing the latter by the NE. The latter leads to an efficient algorithm for calculating the optimal prices.
\end{itemize}
In \cite{candogan2012optimal}, it was shown that the optimal pricing strategy with linear interaction is \emph{independent of the network weights} when the latter is symmetric. With strictly concave interaction, we show that the latter conclusion is not valid. Our results suggest that the price of information (PoI) obeys a nuanced structure in the nonlinear interaction setting.

\section{Network Game \& Optimal Pricing}
\label{sec:model}
Consider a society consisting of $n$ agents, who are connected through a social network represented by the weighted directed graph ${\cal G} = ( {\cal V}, {\cal E} , {\bf G} )$ with ${\cal V} = \{ 1, \ldots, n \}$, ${\cal E} \subseteq {\cal V} \times {\cal V}$ without self-loops. The network is further endowed with a weighted adjacency matrix ${\bf G} \in \RR_+^{n \times n}$. For any  $(j,i) \in {\cal E}$, the $(i,j)$th element of ${\bf G}$ satisfies $G_{ij} > 0$ and the latter quantifies the influence strength of agent $i$ received from agent $j$. If $(j,i) \notin {\cal E}$, then $G_{ij} = 0$. {Note that $G_{ii} = 0$ for $i \in {\cal V}$.}

For each $i \in {\cal V}$, agent $i$ decides on an action $x_i \in \RR_+$. For example, in the context of a consumer social network, this action corresponds to the consumption level of a divisible good. We focus on games with \emph{quadratic payoffs} such that $x_i$ is affected by two sources of externalities. The first one is an agent-specific price, $p_i \in \RR$, decided by a monopolist (can be negative). The second one is due to peer influences where an increase in the neighbors' action leads to a positive effect on $x_i$, i.e., strategic complements \cite{jackson2015games, belhaj2014network}. Formally, the latter effects are described by aggregating the peers' actions exerted through a possibly nonlinear interaction function\footnote{We remark that it is easy to extend our analysis to scenarios when the interaction function is agent dependent.} $f: \RR_+ \to \RR_+$, yielding $\sum_{j = 1}^n G_{ij} f( x_j )$. 

Combining the effects of externalities above leads to the following linear-quadratic payoff in $x_i$: 
\begin{align}
    \label{eq:utility_function}
    u_i(x_i, {\bm x}_{-i} ; p_i) = \left( \sum_{j=1}^{n} G_{ij} f(x_j) + a_i - p_i \right) x_i  - b_i x_i^2.
\end{align}
For agent $i$, $a_i > 0$ {is the marginal incentive parameter} and $b_i > 0$ {is a scaling parameter related to $f(\cdot)$}. We have defined the shorthand notation ${\bm x}_{-i}$ as the vector ${\bm x} := (x_1, \ldots, x_n)$ whose $i$th element $x_i$ is removed. 

In the case when $f(x) = x$, previous works such as \cite{ballester2006s, candogan2012optimal} have studied the linear quadratic game. We focus on the scenario when $f$ is nonlinear (and in particular, strictly concave). The latter models the situation when the peer influences has a diminishing marginal effect which inhibits vigorous actions, see \cite{zhou2016game, cooley2007desegregation,  parise2019variational}.
To this regard, similar payoff functions to \eqref{eq:utility_function} have been studied in \cite{naghizadeh2017on, demange2017optimal} which studied a similar case to ours when the interaction function is applied on the aggregated $x_j$, i.e., $f(\sum_{j=1}^n G_{ij} x_j)$, and in \cite{zhou2016game, parise2019variational} whose general framework includes non linear quadratic games. 

Denote the vectors ${\bm a} = (a_1,\ldots,a_n)$, ${\bm b} = (b_1,\ldots,b_n)$. Given ${\bm p} = (p_1,\ldots,p_n)$, we define a linear-quadratic game by ${\sf Game}({\bm a},{\bm b},{\bf G},f; {\bm p})$ whose payoff function is given in \eqref{eq:utility_function}. We define the \emph{Nash Equilibrium} (NE) of the game as a set of actions, ${\bm x}^{\sf NE} = (x_1^{\sf NE}, \ldots, x_n^{\sf NE} )$, satisfying:
\begin{align}
\label{eq:ne}
	x_i^{\sf NE} \in \argmax_{x_i \in \RR_+ } ~u_i( x_i, {\bm x}_{-i}^{\sf NE} ; p_i),~\forall~i \in {\cal V},
\end{align}
{i.e., when none of the agent has the intention to change his/her action. In Sec.~\ref{sec:NE}, we will show that for any given price vector, the NE is well defined under mild conditions.}

In this paper, we analyze the \emph{optimal pricing} problem for the linear quadratic game \eqref{eq:ne}. The optimal pricing problem involves the monopolist (as {\bf leader}) who maximizes his/her total revenue by setting agent-specific prices, and the agents (as {\bf followers}) whose actions are determined as the NE of ${\sf Game}( {\bm a}, {\bm b}, {\bf G}, f; {\bm p} )$ given the  prices set.
As the total revenue is decided by the actions of agents and the prices set for the agents, this leads to a \emph{bilevel optimization problem}:
\begin{align}
    \max_{ {\bm p}, {\bm x}^{\sf NE} \in \RR^n }~ & {\bm p}^\top {\bm x}^{\sf NE} \label{eq:opt-price} ~~\text{s.t.}~~ 
    {\bm x}^{\sf NE}~\text{satisfies \eqref{eq:ne} with ${\bm p}$,}
    % {\bm x}^{\sf NE}~\text{is NE of}~{\sf Game}( {\bm a}, {\bm b}, {\bf G}, f ; {\bm p} ).
    % \text{satisfies \eqref{eq:ne} given ${\bm p}$} .
    % \notag 
\end{align}
where ${\bm x}^{\sf NE}$ in \eqref{eq:ne} is assumed to be unique.

In the context of consumer social networks, solving \eqref{eq:opt-price} gives a pricing strategy which maximizes the revenue made by the monopolist.
The case of $f(x) = x$ has been studied in \cite{candogan2012optimal}, and the case with nonlinear function applied on aggregated actions is studied in \cite{demange2017optimal} with a different objective of ${\bf 1}^\top {\bm x}^{\sf NE}$. 

\section{Main Results} \label{sec:main}
Our analysis on \eqref{eq:opt-price} is developed via two steps. First, we prove a sufficient condition for the existence of a unique NE. Second, we derive properties on the optimal pricing problem and compare the solution with benchmark strategies. 

\subsection{Nash Equilibrium}
\label{sec:NE}
As a pre-requisite for tackling \eqref{eq:opt-price}, we shall study cases when the NE is well-defined. From \eqref{eq:ne}, the action vector ${\bm x}^{\sf NE}$ is an NE if and only if for any $i \in {\cal V}$,
\begin{align} \label{eq:ne_eqn}
    x_i^{\sf NE} = T_i( {\bm x}^{\sf NE} ) := \max \left\{ 0, \frac{\sum_{j=1}^n G_{ij} f( x_j^{\sf NE} ) + a_i - p_i}{2b_i} \right\} ,
\end{align}
In other words, ${\bm x}^{\sf NE}$ is the fixed point of the nonlinear map ${\bm x} \to T({\bm x}) := [ T_i( {\bm x} ) ]_{i=1}^n$. 
We require the conditions:
\begin{assumption} \label{assu:lips} The following holds:
\begin{enumerate}
\item The interaction function $f$ is $\alpha$-Lipschitz, i.e., for any $x,y \in \RR_+$, it holds $| f(x) - f(y) | \leq \alpha |x-y|$.
\item For any $i \in {\cal V}$, it holds $2 b_i > \alpha \sum_{j=1}^n G_{ij}$.
\end{enumerate}
\end{assumption}
\noindent Note that Assumption~\ref{assu:lips}-1) is satisfied by $f$ which has a bounded derivative over $\RR_+$. On the other hand, Assumption~\ref{assu:lips}-2) requires the curvature parameter of each agent to overcome their own weighted in-degree, scaled by the Lipschitz constant of $f$. Our first result gives a sufficient condition for the existence and uniqueness of NE,
\begin{lemma}
	\label{lem:ne}
	Under Assumption \ref{assu:lips}, Eq.~\eqref{eq:ne_eqn} has a unique fixed point. As such, ${\sf Game}({\bm a},{\bm b},{\bm G}, f; {\bm p} )$ has a unique NE.
\end{lemma}
\noindent 
See Appendix~\ref{app:ne} for the proof.
% The proof is relegated to Appendix~\ref{app:ne}. 
When $f(x) = x$, we have $\alpha=1$ and our result recovers \cite[Theorem 1]{candogan2012optimal}. The key difference lies in the case when $\alpha \neq 1$, where the condition in Assumption~\ref{assu:lips}-2) will be affected by the maximum derivative of $f$.

Observe that ${\sf Game}( {\bm a}, {\bm b}, {\bf G}, f ; {\bm p} )$ is subject to scaling ambiguity. E.g., the interaction function can be shifted or scaled arbitrarily. To this regard, it is better to work with a `normalized' instance of the linear quadratic game. Two games are \emph{equivalent} if their NEs coincide.
{With a direct substitution of variables, we} observe that: 
\begin{corollary} \label{lem:normalize}
For any game ${\sf Game}({\bm a},{\bm b}, {\bf G},f ; {\bm p})$ that satisfies Assumption \ref{assu:lips}, there exists an equivalent game ${\sf Game}({\hat{\bm a}},{\hat{\bm b}},\hat{\bf G},  \hat{f}; {\bm p} )$ with the parameters
\[
\hat{\bm a} = {\bm a} + f(0) {\bf G} {\bf 1},~\hat{\bm b} = {\bm b},~\hat{f}(x) = \frac{ f(x) - f(0) }{\alpha},~\hat{\bf G} = \alpha {\bf G}.
\]
Moreover, we have $\hat{f}(0) = 0, \hat{f}(x) \leq x$ for any $x \in \RR_+$.
\end{corollary}

\noindent For the rest of this paper, we will study the linear quadratic games with a normalized interaction function such that $f(x)$ is $1$-Lischitz and satisfies $f(0)=0$, $f(x) \leq x$ for any $x \in \RR_+$. This is without loss of generality due to Corollary~\ref{lem:normalize}. As a concrete example, one may consider $f(x) = \ln(1+x)$. From now on, we shall denote ${\bm x}^{\sf NE}( {\bm p} )$ as the NE action vector to ${\sf Game}( {\bm a}, {\bm b}, {\bf G}, f ; {\bm p} )$ with the normalized interaction function.

\subsection{Optimal Pricing} 
We present algorithmic and analytical results on the bilevel optimization problem \eqref{eq:opt-price}. Define the diagonal matrix $\bm{\mathcal{B}} := 2{\rm Diag}( b_1, \ldots, b_n )$ and consider the following set of assumptions that have been strengthened from Assumption~\ref{assu:lips}:
\begin{assumption} \label{assu:f_strong}
    The interaction function $f$ is $1$-Lipschitz, non-decreasing, twice differentiable, concave on $\RR^+$ and satisfies $f(0) = 0$, $f(x) \leq x$ for any $x \in \RR^+$.
\end{assumption}
\begin{assumption} \label{assu:B}
There exists $\rho, \rho' > 0$ such that the matrices 
$\bm{\mathcal{B}} - \frac{1}{2} ( {\bf G} + {\bf G}^\top) - \rho {\bf I}$, $\bm{\mathcal{B}} - {\bf G} - \rho' {\bf I}$ are diagonally dominant\footnote{${\bf A} \in \RR^{n \times n}$ is diagonally dominant if $|A_{ii}| \geq \sum_{j \neq i} |A_{ij}|$ for all $i$.}.
\end{assumption}

\subsubsection{Efficient Computation} To derive an efficient algorithm for solving \eqref{eq:opt-price}, our approach consists of showing that the optimal pricing problem is equivalent to a strongly concave optimization problem on the action vector ${\bm x}$.

Our approach hinges on establishing the correspondence between an NE and a price vector. However, the nonlinearity of the $\max\{0, \cdot\}$ map in \eqref{eq:ne_eqn} could hinder our development. Fortunately, we show that at the optimal pricing strategy, the effect of the $\max\{0, \cdot\}$ map can be ignored: 
\begin{lemma} \label{lem:positive_x}
Let ${\bm p}^\star$ be an optimal solution to \eqref{eq:opt-price}. Then {the NE to ${\sf Game}( {\bm a}, {\bm b}, {\bf G}, f ; {\bm p}^\star )$ must satisfy} ${\bm x}^{\sf NE} ({\bm p}^\star ) > {\bm 0}$.
\end{lemma}
\noindent 
{The proof is modified from \cite[Lemma 7]{candogan2012optimal} by adapting the latter to $f(\cdot)$.}
Together with \eqref{eq:ne_eqn}, Lemma~\ref{lem:positive_x} implies that at the optimal price ${\bm p}^\star$, the corresponding NE must satisfy: 
\begin{align} \label{eq:ne_eqn_opt}
    {\bm x}^{\sf NE} ({\bm p}^\star) = \bm{\mathcal{B}}^{-1} \big( {\bf G} {\bf f}( {\bm x}^{\sf NE} ({\bm p}^\star) ) + {\bm a} - {\bm p}^\star \big) > {\bm 0} ,
\end{align}
where we have defined ${\bf f}( {\bm x} ) = ( f( x_1 ) , \ldots, f( x_n ) ) \in \RR^n$.
This suggests us to rewrite {the objective function ${\bm p}^\top {\bm x}$ in the bilevel optimization problem \eqref{eq:opt-price} by substituting ${\bm p} = {\bm a} + {\bf G} {\bf f}( {\bm x}) - \bm{\mathcal{B}} {\bm x}$. We consider}:
\begin{align} \label{eq:jx}
    J( {\bm x} ) := {\bm x}^\top \left( {\bm a} + {\bf G} {\bf f}( {\bm x} ) - \bm{\mathcal{B}} {\bm x} \right) ,
\end{align}
By Lemma~\ref{lem:positive_x}, $J( {\bm x}^{\sf NE} ({\bm p}^\star) )$ is the optimal total revenue achieved by solving \eqref{eq:opt-price}. We further note that
\begin{proposition} \label{prop:strong}
    Under Assumptions~\ref{assu:f_strong} and \ref{assu:B}, the function $J({\bm x})$ is strongly concave with modulus $2\rho$, and {$ {\bm x}^{\sf NE} ({\bm p}^\star) $ satisfying \eqref{eq:ne_eqn_opt} is the unique maximizer of $J( {\bm x} )$ over ${\bm x} \in \RR_{+}^n$}.
\end{proposition}

Lemma~\ref{lem:positive_x} and Proposition~\ref{prop:strong} imply that \eqref{eq:opt-price} is equivalent to a \emph{single-level} optimization problem $\max_{ {\bm x} \in \RR^n_+ } J({\bm x})$.
A standard procedure is to consider a projected gradient (PG) method for solving \eqref{eq:opt-price}: for $t = 0, 1, \ldots$
\begin{align} \label{eq:pg}
    {\bm x}^{t+1} = {\max} \big\{ {\bm 0}, {\bm x}^t + \gamma \nabla J( {\bm x}^t ) \big\},
\end{align}
where $\gamma > 0$ is a step size. Assume in addition that $f''(x) \geq -M$, it can be shown that $\grd J({\bm x})$ is $\nu$-Lipschitz, see \eqref{eq:bound_nablaJ}. Setting $\gamma = 1/\nu$ guarantees that \eqref{eq:pg} converge to the optimal action vector ${\bm x}^\star$ geometrically \cite[Theorem 10.29]{beck2017first}.
The optimal prices can be found by ${\bm p}^\star = {\bm a} + {\bf G} {\bf f}( {\bm x}^\star ) - \bm{\mathcal{B}} {\bm x}^\star$.

\subsubsection{Network Externality \& Price Discrimination} \label{sec:opt-price} To simplify notation, let ${\bm p}^\star$ be an optimal pricing vector and ${\bm x}^\star := {\bm x}^{\sf NE}( {\bm p}^\star )$ be the NE action vector. We aim to compare ${\bm p}^\star$ to a
\emph{network agnostic} strategy where network externalities are ignored by the monopolist. In this case, the monopolist models the payoff function as \eqref{eq:utility_function} without $\sum_{j = 1}^n G_{ij} f(x_j)$. The optimal pricing problem will be reduced to 
$\max_{ {\bm p} \in \RR^n } \sum_{i=1}^n \frac{p_i}{2b_i} \max \{ 0, a_i - p_i \}$,
whose optimal solution is given by ${\bm p}_0 := {\bm a}/2$. We further denote ${\bm x}_0 = {\bm x}^{\sf NE} ( {\bm p}_0 )$ as the corresponding NE.
To this end, we are interested in the \emph{price of information (PoI)}:
\beq 
PoI := J( {\bm x}^\star ) / J( {\bm x}_0 ) ,
\eeq 
which is the ratio between the optimal revenue and the revenue with network agnostic strategy. 

The first result is to show that in general, the above \emph{network agnostic} pricing strategy is strictly \emph{suboptimal}:
\begin{theorem}\label{the:the2} Under Assumptions \ref{assu:f_strong} and \ref{assu:B}. If $f$ is strictly concave, then it holds
\beq \label{eq:strict_bound}
J({\bm x}^\star) = ( {\bm x}^\star )^\top {\bm p^\star} > ( {\bm x}_0 )^\top {\bm p}_0 = J( {\bm x}_0 ).
\eeq 
Furthermore, if there exists $M \geq 0$ such that $f''(x) \geq -M$ for any $x \geq 0$, then the PoI satisfies
\beq \label{eq:tight_bound}
\frac{J( {\bm x}^\star ) }{J( {\bm x}_0 ) } -1 \geq \frac{4 \rho}{ \nu^2 } \frac{ \left( {\bm a}^\top \bm{\mathcal{B}}^{-1} {\bf G} \, {\bf h}( \frac{1}{2} \bm{\mathcal{B}}^{-1} {\bm a} ) \right)^2} { ( {\bm a}^\top ( \bm{\mathcal{B}} - {\bf G} )^{-1} {\bm a}) \| ( \bm{\mathcal{B}} - {\bf G} )^{-1} {\bm a} \|^2 },
\eeq 
where ${\bf h}( {\bm x} ) := {\bf f}( {\bm x} ) - D_f( {\bm x} ) {\bm x}$ with $D_f( {\bm x} ) := {\rm Diag}( f'(x_1), \ldots, f'(x_n) )$, and 
\begin{align} \label{eq:nu_def}
\nu := \| 4 {\bm b} + ( {\bf G} + {\bf G}^\top) {\bf 1} + M {\bf G}^\top \bar{\bm x}^{\max} \|_\infty,
\end{align}
where $\bar{\bm x}^{\max} = \max\{ \frac{1}{2} ( \bm{\mathcal{B}} - {\bf G} )^{-1} {\bm a}, \frac{1}{2} ( \bm{\mathcal{B}} - \frac{ {\bf G} + {\bf G}^\top }{2} )^{-1} {\bm a} \}$ is taken with the element-wise maximum.
\end{theorem}
\noindent 
{The proof is established by applying Taylor’s theorem on
$J( {\bm x}_0 )$ and relate it to $\| {\bm x}^\star - {\bm x}_0\|$;} see Appendix~\ref{app:thm2}. 
Our result shows that for \emph{strictly concave interaction function}, it is necessary for the monopolist to explore \emph{network knowledge} to find the maximum total revenue through \emph{pricing discrimination}. 

We first compare our result to that of \cite{candogan2012optimal} which studied the case of $f(x) = x$, i.e., the interaction function is linear. Corollary 1 therein shows that when ${\bf G} = {\bf G}^\top$, then the \emph{network agnostic} pricing strategy ${\bm p}_0 = {\bm a}/2$ is \emph{optimal}. In the special case with uniform incentive ${\bm a} = a {\bf 1}$, there is no need for pricing discrimination for symmetric network. Note that this is compatible with our bound on PoI in \eqref{eq:tight_bound} as ${\bf h}({\bm x}) = {\bm 0}$ when $f(x) = x$. On the other hand, we have ${\bf h}({\bm x}) > {\bm 0}$ when $f$ is strictly concave. For the latter case, the lower bound in \eqref{eq:tight_bound} will be strictly positive. This implies that \emph{regardless of whether ${\bf G} = {\bf G}^\top$ or not}, it is strictly beneficial for the monopolist to explore a network externality dependent pricing strategy.

Furthermore, the lower bound for PoI in \eqref{eq:tight_bound} quantifies the minimum gain in total revenue with the optimal pricing vector ${\bm p}^\star$ in the presence of network information. We first notice that this bound (i) increases with $\rho$ in Assumption~\ref{assu:B}, and (ii) decreases with the curvature of $f$ through $M$. Furthermore, from the right hand side of \eqref{eq:tight_bound}, the bound roughly behaves as a concave function of ${\bf G}^\top \bm{\mathcal{B}}^{-1} {\bm a}$, suggesting a complex interaction between the incentive level ($\bm{\mathcal{B}}^{-1} {\bm a}$) and the network structure (${\bf G}$). We will further study this behavior in Sec.~\ref{sec:exp} using numerical experiments.

The second scenario is to consider a case when the monopolist may process the knowledge of the network structure but is unable to offer price discrimination due to the cost of implementation. We aim to quantify the value of network structure even without the ability to charge prices differently for different agents. To gain insights, we focus on a `symmetric' game ${\sf Game}({\bm a}, {\bm b}, {\bf G}, f; {\bm p})$ satisfying 
\begin{align} \label{eq:unif_game}
{\bm a} = \bar{a} {\bf 1},~{\bm b} = \bar{b} {\bf 1},~{\bf G} {\bf 1} = \bar{g} {\bf 1}.
\end{align}
% The following is a consequence of the uniqueness of the NE:
\begin{corollary} \label{cor:symmetric}
Under Assumptions~\ref{assu:f_strong}, \ref{assu:B}. Suppose that the game satisfies \eqref{eq:unif_game}. Given a uniform price strategy $\bar{p} {\bf 1}$, the NE action vector must be uniform, i.e., ${\bm x}^{\sf NE}( \bar{p} {\bf 1} ) = \bar{x} {\bf 1}$, $\bar{x} \geq 0$.
\end{corollary}
\noindent 
{The corollary is obtained by showing that $\bar{x} {\bf 1}$ is an NE, followed by its uniqueness property.}
The corollary implies that when the optimal pricing problem is restricted to finding a uniform price, the bilevel optimization problem is reduced into a \emph{one dimensional optimization} $\max_{ \bar{x} \geq 0 } J( \bar{x} {\bf 1} )$ [cf.~\eqref{eq:jx}].

Let $\bar{p}^\star$ to be the \emph{optimal uniform price} and $\bar{x}^\star {\bf 1}$ be the uniform NE  action which satisfies $\bar{x}^\star = \argmax_{ \bar{x} > 0 } J( \bar{x} {\bf 1} )$. Similar to the previous scenario, the uniform network agnostic optimal price is given by $\bar{p}_0 = \bar{a}/2$. Denote $\bar{x}_0 {\bf 1}$ as the uniform NE action for the latter. {Using a similar proof idea as Theorem~\ref{the:the2},} we bound the gain of $\bar{p}^\star$ over $\bar{p}_0$:
\begin{theorem}
    \label{the:the3}
    Under Assumptions \ref{assu:f_strong}, \ref{assu:B}. Suppose that the game satisfies \eqref{eq:unif_game}. Let $h(x) := { f(x) - f'(x) x }$, it holds 
    \begin{align} \label{eq:tight_bound_scalar}
    \frac{J( \bar{x}^\star {\bf 1} )}{ J(\bar{x}_0 {\bf 1} )} - 1 \geq \left(  \frac{ \bar{g} }{ \bar{a} }  \left(1 - \frac{ \bar{g} }{ 2 \bar{b} } \right)  h\left(\frac{ \bar{a} }{ 4 \bar{b} } \right) \right)^2.
    \end{align}
\end{theorem}

Obtaining $\bar{x}^\star$ and thus $\bar{p}^\star$ requires the knowledge of the network externality parameters $\bar{g}$, $\bar{a}$, $\bar{b}$, $f(\cdot)$. Notice that the detailed structure of ${\bf G}$ is not required. Meanwhile, the benchmark price $\bar{p}_0$ requires only $\bar{a}$.
For strictly concave $f(\cdot)$, we have $h(x) > 0$ for any $x > 0$. As such, similar to \eqref{eq:tight_bound}, the above continues to show that the \emph{optimal uniform price} informed by network externalities is a strictly better pricing strategy than the network agnostic strategy, regardless of whether ${\bf G}$ is symmetric or not. Finally, we note that the left hand side of \eqref{eq:tight_bound_scalar} is a lower bound to the PoI.

\section{Numerical Experiment} \label{sec:exp}
This section presents numerical experiments on how network information will influence the total revenue in optimal pricing and validate the results in Sec.~\ref{sec:main}.

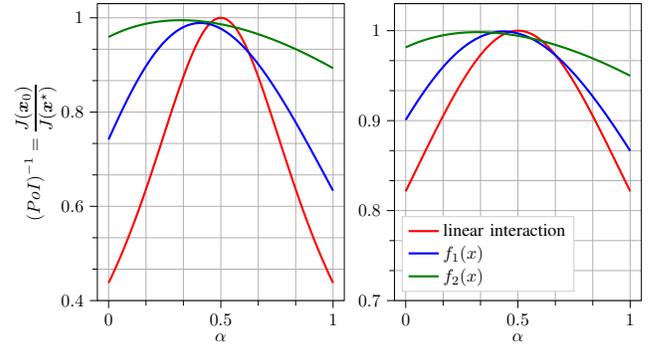
\begin{figure}
\centering
\resizebox{0.95\linewidth}{!}{\input{./figs_tex/star_pa}}
\vspace{-0.3cm}
 \caption{(Left) Star graph with $a = 1; b = 1; n = 10$. Maximum PoI achieved by $f_1(\cdot)$ is $\frac{1}{0.989}$ at $\alpha = 0.41$, $f_2(\cdot)$ is $\frac{1}{0.995}$ at $\alpha = 0.32$; (Right) PA graph with $a = 1; b = 2; n = 100$. Maximum PoI achieved by $f_1(\cdot)$ is $\frac{1}{0.999}$ at $\alpha = 0.44$, $f_2(\cdot)$ is $\frac{1}{0.998}$ at $\alpha = 0.32$.}
\vspace{-.3cm}
\label{fig:star_pa}
\end{figure}

\subsection{Network Symmetricity}
Our first set of experiments examines the ratio $\frac{ J( {\bm x}_0 ) }{ J( {\bm x}^\star ) }$, i.e., the total revenues under network agnostic pricing strategy ${\bm p}_0$ [cf.~Sec.~\ref{sec:opt-price}] and the optimal pricing strategy ${\bm p}^\star$ solving \eqref{eq:opt-price}.
We compare this ratio of revenues when different interaction functions $f(\cdot)$ are used in \eqref{eq:utility_function}: the linear function $f(x)=x$, or the normalized strictly concave functions $f_1(x) = \ln(x+1)$, $f_2(x) = 0.1\ln(1+10x)$.\vspace{.1cm}

\noindent \textbf{Star Graph}: Consider a star graph topology described by the asymmetric matrix ${\bf G}^{\sf star}$ with $G_{ij}^{\sf star} = 1$ if $i=1, j \neq i$, otherwise $G_{ij}^{\sf star} = 0$ such that agent $1$ is the `central' agent in the network. For the network game, we use the weighted adjacency matrix ${\bf G} = \alpha {\bf G}^{\sf star} + (1-\alpha) ({\bf G}^{\sf star})^\top$ parameterized by $\alpha \in [0,1]$. Increasing $\alpha$ trades off the impact level of the central agent has on her neighbors and the impact level of her neighbors on the central agent. At $\alpha = 0.5$, the directed impact levels are balanced and ${\bf G} = {\bf G}^\top$ is a symmetric adjacency matrix. 
With $n=10$ agents, Fig.~\ref{fig:star_pa} (left) shows the inverse PoI against $\alpha \in [0,1]$. We first notice that unlike the case of linear interaction which achieves $\frac{ J( {\bm x}_0 ) }{ J( {\bm x}^\star ) } = 1$ at $\alpha = 0.5$, the two cases of nonlinear interaction functions are skewed where the maximum ratios of revenues are achieved at $\alpha < 0.5$, and the ratios are strictly less than $1$. These corroborate Theorem~\ref{the:the2}. 
% Importantly, in extreme "asymmetric" scenarios (\(\alpha = 0\), \(\alpha = 1\)), the availability of network information enhances the total revenue noticeably.
\vspace{.1cm}

\noindent \textbf{Preferential Attachment (PA) Graph}:
We consider the network specified by a PA graph with $n=100$ agents. We first generate an upper triangular matrix ${\bf G}^{\sf pa}$ using {the procedure in \cite[Section 5]{candogan2012optimal}}. Similar to the previous example, we use the weighted adjacency matrix ${\bf G} = \alpha {\bf G}^{\textsf{pa}} + (1 - \alpha)({\bf G}^{\textsf{pa}})^\top$. 
Similar to our findings in the star graph example, Fig.~\ref{fig:star_pa} also shows that the maximum inverse PoI of revenues is achieved at $\alpha < 0.5$ for nonlinear interaction.

{For both examples, the PoIs under concave interaction function is reduced at the extreme cases ($\alpha \to 0$ or $\alpha \to 1$) when compared to when $f(x) = x$. We note that this is supported by a related observation in \cite{demange2017optimal} that shows the targeting strategy is less effective with concave interaction.}

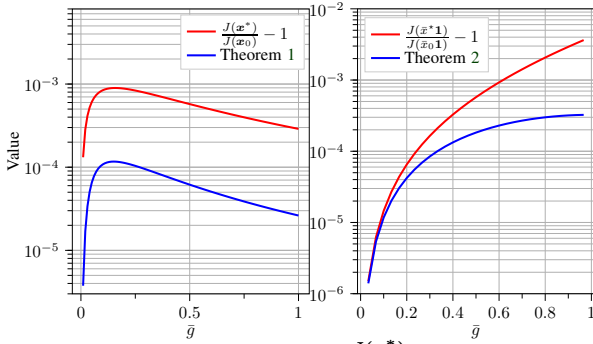
\begin{figure}
\centering
\resizebox{0.9\linewidth}{!}{\input{./figs_tex/bound_all}}
\vspace{-0.3cm}
\caption{(Left) Comparing the actual ratio $\frac{J({\bm x}^*)}{J({\bm x}_0)}-1$ to the lower bound in Theorem~\ref{the:the2} with $a=2,b =\bar{g}$; (Right) Comparing the ratio with uniform prices $\frac{J(\bar{x}^\star {\bm 1} )}{J(\bar{x}_0 {\bm 1} )}-1$ to the lower bound in Theorem~\ref{the:the3} with $a=2,b =1$.}
\vspace{-.4cm}
\label{fig:all_bound}
\end{figure}

\subsection{Lower Bounds of Network Information Gap}
Our second set of experiments evaluates the tightness of the lower bounds established in Theorems~\ref{the:the2}, \ref{the:the3}. We focus on the case with the interaction function $f_1(x) = \ln(1+x)$. With the parameter $\bar{g} \geq 0$, the social network is described by a ring graph such that $G_{ij}^{\text{ring}} = \bar{g}$ if $i \equiv (j+1)~{\rm mod}~n$; otherwise $G_{ij}^{\text{ring}} = 0$. We consider adjusting the interaction strength $\bar{g}$ from $0$ to $1$ and compare $\frac{ J( {\bm x}^\star ) }{ J( {\bm x}_0 ) } - 1$ to the lower bounds proven in Theorems~\ref{the:the2}, \ref{the:the3}. We also fix the number of agents at $n=100$.

Fig.~\ref{fig:all_bound} (left) compares the bound in Theorem~\ref{the:the2} under the general setting where price discrimination is allowed. We observe that our theoretical bound captures the concave nature of the gap $\frac{ J( {\bm x}^\star ) }{ J( {\bm x}_0 ) } - 1$ w.r.t.~$\bar{g}$. However, we also observe that the bound is in general loose. Similarly, Fig.~\ref{fig:all_bound} (right) compares the bound in Theorem~\ref{the:the3} which is specialized to the uniform pricing strategy setting. Here, we compare $\frac{J(\bar{x}^\star \bm{1} )}{J(\bar{x}_0 \bm{1} )}-1$ to the theoretical bound in the theorem. The figure shows that the proposed bound is tighter in this case.\vspace{-.1cm}

\vspace{-.1cm}
\section{Conclusions}
\label{sec:concl}
We considered the optimal pricing problem for linear quadratic game with externalities exerted through strictly concave interaction. 
% We show that the optimal pricing problem can be efficiently solved as a strongly concave maximization. Moreover, w
We showed that it is \emph{necessary} for the monopolist to use knowledge on network structure in finding the optimal pricing strategy. It contrasts with the \emph{linear interaction} case, where the network structure knowledge can be \emph{unnecessary} in symmetric network.
Future works include extension to general games with other forms of nonlinearity.\vspace{-.1cm}

\appendix \vspace{-.1cm}

\allowdisplaybreaks

\subsection{Proof of Lemma~\ref{lem:ne}} \label{app:ne}
\begin{proof} 
We have\vspace{-.2cm}
\begin{align*}
& \left\| T( {\bm x} ) - T( {\bm y} ) \right\|^2 \leq \sum_{i = 1}^n \left( \sum_{j = 1}^n \frac{G_{ij} (f(x_j)-f(y_j))}{2b_i} \right)^2 \\
& \leq \alpha^2 \sum_{i = 1}^n \left( \sum_{j = 1}^n \frac{ G_{ij} |x_j - y_j| }{2b_i} \right)^2 \leq \alpha^2 \big( \rho(\bm{\mathcal{B}}^{-1} {\bf G}) \, \| x - y \| \big)^2,
\end{align*}
where $\rho(\cdot)$ denotes the spectral radius for a matrix.
With a slight modification to \cite[Lemma 4]{candogan2012optimal}, it can be shown that $\alpha \, \rho(\bm{\mathcal{B}}^{-1} {\bf G}) < 1$. {Together, this shows that $T(\cdot)$ is a contractive map. Applying the Banach's fixed point theorem yields the lemma.}
\end{proof}

\subsection{Proof of Proposition~\ref{prop:strong}} \label{app:strong}
\begin{proof} 
To prove the first part of the proposition, we notice 
\begin{align} \label{eq:grad_J}
\nabla J({\bm x} ) & = {\bm a} - 2 \bm{\mathcal{B}} {\bm x} + {\bf G} {\bf f} ({\bm x}) + D_f({\bm x}){\bf G^\top} {\bm x},
\end{align}
where $D_f({\bm x}) := {\rm Diag} ( f'(x_1), \dots, f'(x_n) )$. Furthermore,
\begin{equation}
\begin{aligned}
\textstyle \frac{\partial^2 J}{\partial x_i^2}   & \textstyle = - 4b_i + \sum_{j} G_{ji} x_j f''(x_i) \leq -4b_i \\
\textstyle \frac{\partial^2 J}{\partial x_i x_j} &= G_{ij} f'(x_j) + G_{ji} f'(x_i),~i \neq j.
\end{aligned}
\end{equation}
Observe that as $f$ is strictly concave, 
\begin{align}
\textstyle \left| \frac{\partial^2 J}{\partial x_i^2} \right| \geq 4b_i - \sum_{j} G_{ji} x_j f''(x_i) \geq 4b_i.
\end{align}
Furthermore, as $f(x)$ is $1$-Lipschitz, we have $|f'(x)| \leq 1$. By Assumption~\ref{assu:B}, we have
\begin{equation} \notag
\begin{aligned}
    & \textstyle 4b_i - 2\rho \geq \sum_{j=1}^n ( G_{ij} + G_{ji} ) \\
    & \textstyle \geq \sum_{j \neq i} \big| G_{ij} f'(x_j) + G_{ji} f'(x_i) \big| = \sum_{j \neq i} \left| \frac{\partial^2 J}{\partial x_i x_j} \right|.
\end{aligned}
\end{equation}
As such, $\nabla^2 J( {\bm x} ) + 2\rho {\bf I}$ is diagonally dominant. 
By the Gershgorin Circle Theorem, $\nabla^2 J( {\bm x} ) + 2\rho {\bf I}$ is negative semidefinite. We conclude that $J( {\bm x} )$ is $2\rho$-strongly concave. 

For the second part of the proposition, by \eqref{eq:grad_J}, the unique maximizer, ${\bm x}^\star$, of $J({\bm x})$ must be strictly positive. Since ${\bm x} > 0$, this must correspond to a pricing vector ${\bm p} = {\bm a} + {\bf G} f( {\bm x}^\star ) - \bm{\mathcal{B}} {\bm x}^\star$  feasible to \eqref{eq:opt-price}. Then, such pricing vector must also be optimal to \eqref{eq:opt-price} and show that $J( {\bm x}^\star ) = J( {\bm x}^{\sf NE}( {\bm p}^\star ) )$. 
\end{proof}

\subsection{Proof of Theorem~\ref{the:the2}}\label{app:thm2}
\begin{proof}
Note that at the pricing strategy of ${\bm p}_0 = {\bm a}/2$, one must have ${\bm x}_0 := {\bm x}^{\sf NE}( {\bm p}_0 ) > {\bm 0}$ and thus 
\begin{align} \label{eq:x0_first}
    {\bm x}_0 = \bm{\mathcal{B}}^{-1} \left( \frac{\bm a}{2} + {\bf G} {\bf f}( {\bm x}_0 ) \right).
\end{align}
Denoting ${\bf h}( {\bm x} ) := {\bf f}( {\bm x} ) - D_f( {\bm x} ) {\bm x}$, where $D_f( {\bm x} ) := {\rm Diag}( f'(x_1), \ldots, f'(x_n) )$. We observe that ${\bf h}( {\bm x} ) \geq {\bm 0}$ for any ${\bm x} > {\bm 0}$, and it is an \emph{non-decreasing} function in ${\bm x}$ due to the concavity of $f$. We have 
\begin{equation}
\begin{aligned}
    \nabla J( {\bm x}_0 ) & = {\bm a} - 2 \bm{\mathcal{B}} {\bm x}_0 + {\bf G} {\bf f} ({\bm x}_0) + D_f({\bm x}_0){\bf G^\top} {\bm x}_0 \\
    & = - {\bf G} {\bf f}( {\bm x}_0 ) + D_f({\bm x}_0){\bf G^\top} {\bm x}_0.
\end{aligned}
\end{equation}
Furthermore, if $f$ is strictly concave, then ${\bf h}( {\bm x} ) > {\bm 0}$ and we further note that
\begin{equation} \label{eq:x0nablax0}
\begin{aligned}
{\bm x}_0^\top \nabla J( {\bm x}_0 ) & = - {\bm x}_0^\top {\bf G} {\bf h}({\bm x}_0) < 0.
\end{aligned}
\end{equation}
The above establishes $\nabla J( {\bm x}_0 ) \neq {\bm 0}$. Since we also have ${\bm x}_0 > {\bm 0}$, this shows the suboptimality of ${\bm x}_0$ and proves \eqref{eq:strict_bound}.

To estimate the ratio $J( {\bm x}^\star ) / J( {\bm x}_0 )$, our idea is to apply Taylor's theorem to approximate $J( {\bm x}_0 )$ around ${\bm x}^\star$. To handle the nonlinearity, we shall first derive lower and upper bound estimates for ${\bm x}_0$, ${\bm x}^\star$ and exploit the monotonicity of $f(\cdot)$. 

From \eqref{eq:x0_first}, we immediately obtain ${\bm x}_0 > \frac{1}{2} \bm{\mathcal{B}}^{-1} {\bm a}$.
% \begin{equation}
% {\bm x}_0 = \bm{\mathcal{B}}^{-1} \left( \frac{\bm a}{2} + {\bf G} {\bf f}( {\bm x}_0 ) \right) > \frac{1}{2} \bm{\mathcal{B}}^{-1} {\bm a}
% \end{equation}
On the other hand, using $f(x) \leq x$, we have
\begin{equation}
2 \bm{\mathcal{B}} {\bm x}_0 = {\bm a} + 2 {\bf G} {\bf f}( {\bm x}_0 ) \leq {\bm a} + 2 {\bf G} {\bm x}_0.
\end{equation}
Under Assumption~\ref{assu:B}, $( \bm{\mathcal{B}} - {\bf G} )^{-1}$ is  non-negative. Thus,
\begin{equation} \label{eq:ub_0}
{\bm x}_0 \leq (1/2) ( \bm{\mathcal{B}} - {\bf G} )^{-1} {\bm a}.
\end{equation}
Similarly, we can lower and upper bound ${\bm x}^\star$. From the optimality of ${\bm x}^\star$, we have 
\begin{equation}
    2 \bm{\mathcal{B}} {\bm x}^\star = {\bm a} + {\bf G} {\bf f}( {\bm x}^\star ) + D_f( {\bm x}^\star ) {\bf G}^\top {\bm x}^\star
\end{equation}
We immediately get ${\bm x}^\star \geq \frac{1}{2} \bm{\mathcal{B}}^{-1} {\bm a}$. On the other hand, using $|f'(x)| \leq 1$ and $f(x) \leq x$ yields
\begin{equation}
    2 \bm{\mathcal{B}} {\bm x}^\star \leq {\bm a} + ( {\bf G} + {\bf G}^\top ) {\bm x}^\star
\end{equation}
The matrix $( \bm{\mathcal{B}} - \frac{ {\bf G} + {\bf G}^\top }{2} )^{-1}$ is non-negative. This implies 
\begin{equation} \label{eq:ub_star}
\textstyle {\bm x}^\star \leq \frac{1}{2} \left( \bm{\mathcal{B}} - \frac{ {\bf G} + {\bf G}^\top }{2} \right)^{-1} {\bm a}.
\end{equation}

Our next step is to derive a few estimates using the mean value theorem. To this end, define ${\cal X} = \{ {\bm x}_0 + \alpha( {\bm x}^\star - {\bm x}_0 ) : \alpha \in [0,1] \}$ as the line segment between ${\bm x}_0$ and ${\bm x}^\star$.  
By $|f''(x)| \leq M$, applying the Gershgorin circle theorem shows that for the negative definite matrix $\grd^2 J({\bm x})$, we can lower bound its minimum eigenvalue as
\begin{equation} \label{eq:bound_nablaJ}
    \begin{aligned}
    & \lambda_{\min}( \grd^2 J( {\bm x} ) ) \\
    & \geq - \max_i \left\{ M \sum_{j=1}^n G_{ji} x_j + \sum_{j=1}^n ( G_{ij} + G_{ji} ) + 4b_i \right\} \\
    & = - \| 4 {\bm b} + ( {\bf G} + {\bf G}^\top) {\bf 1} + M {\bf G}^\top {\bm x} \|_\infty = - \upsilon.
    \end{aligned}
\end{equation}
The lower bound in the above decreases when ${\bm x}$ increases. Thus, we conclude that $\lambda_{\min}( \grd^2 J( {\bm x} ) ) \geq - \nu$ for any ${\bm x} \in {\cal X}$, where we recall $\nu$ from \eqref{eq:nu_def}
% \begin{align}
% \nu := \| 4 {\bm b} + ( {\bf G} + {\bf G}^\top) {\bf 1} + M {\bf G}^\top \bar{\bm x} \|_\infty,
% \end{align}
and $\bar{\bm x}^{\max} = \max\{ \frac{1}{2} ( \bm{\mathcal{B}} - {\bf G} )^{-1} {\bm a}, \frac{1}{2} ( \bm{\mathcal{B}} - \frac{ {\bf G} + {\bf G}^\top }{2} )^{-1} {\bm a} \}$ is the element-wise maximum between the upper bounds on ${\bm x}^\star$, ${\bm x}_0$ in \eqref{eq:ub_0}, \eqref{eq:ub_star}.
As a consequence of the above analysis, we note that by the mean value theorem, for some ${\bm x} \in {\cal X}$,
\begin{align}
\grd J( {\bm x}_0 ) = \grd^2 J( {\bm x} ) ( {\bm x}_0 - {\bm x}^\star ).
\end{align}
since $\grd J( {\bm x}^\star ) = {\bm 0}$. It follows that 
\begin{align}
    \| \grd J( {\bm x}_0 ) \| \leq \| \grd^2 J( {\bm x} ) \| \| {\bm x}_0 - {\bm x}^\star \| \leq \nu \| {\bm x}_0 - {\bm x}^\star \|.
\end{align}
On the other hand, using \eqref{eq:x0nablax0}, we have 
\begin{align}
    \| \grd J( {\bm x}_0 ) \| \geq ( \| {\bm x}_0 \|  )^{-1} {\bm x}_0^\top {\bf G} {\bf h}({\bm x}_0).
\end{align}

We conclude the proof by applying Taylor's theorem on $J( {\bm x}_0 )$. For some ${\bm x} \in {\cal X}$, it holds
\begin{equation}
    \begin{aligned}
    J( {\bm x}_0 ) & = J( {\bm x}^\star ) + \frac{1}{2} ( {\bm x}_0 - {\bm x}^\star ) \grd^2 J( {\bm x} ) ( {\bm x}_0 - {\bm x}^\star ) \\
    & \geq J( {\bm x}^\star ) - \rho \| {\bm x}_0 - {\bm x}^\star \|^2.
    \end{aligned}
\end{equation}
Subsequently, 
\begin{equation}
    \begin{aligned}
    \frac{ J( {\bm x}^\star ) }{ J ( {\bm x}_0 ) } - 1 & \geq \rho \frac{ \| {\bm x}_0 - {\bm x}^\star \|^2 }{ J( {\bm x}_0 ) } \geq \frac{\rho}{\upsilon^2} \frac{ \| \grd J( {\bm x}_0 ) \|^2 }{ J( {\bm x}_0 ) }.
    \end{aligned}
\end{equation}
Note that with the upper bound on ${\bm x}_0$, we get the estimates: $J( {\bm x}_0 ) \leq \frac{1}{4} {\bm a}^\top ( \bm{\mathcal{B}} - {\bf G} )^{-1} {\bm a}$, $\| {\bm x}_0 \| \leq \| \frac{1}{2} ( \bm{\mathcal{B}} - {\bf G} )^{-1} {\bm a} \|$, and 
\begin{align}
{\bm x}_0^\top {\bf G} {\bf h}({\bm x}_0) \geq \frac{1}{2} {\bm a}^\top \bm{\mathcal{B}}^{-1} {\bf G} \, {\bf h} \big( {\textstyle \frac{1}{2} \bm{\mathcal{B}}^{-1} {\bm a}} \big).
\end{align}
This yields 
\begin{align}
\frac{ J( {\bm x}^\star ) }{ J ( {\bm x}_0 ) } - 1 \geq \frac{4 \rho}{\nu^2} \frac{ \big( {\bm a}^\top \bm{\mathcal{B}}^{-1} {\bf G} \, {\bf h} ( {\textstyle \frac{1}{2} \bm{\mathcal{B}}^{-1} {\bm a}} ) \big)^2 }{ ( {\bm a}^\top ( \bm{\mathcal{B}} - {\bf G} )^{-1} {\bm a}) \| ( \bm{\mathcal{B}} - {\bf G} )^{-1} {\bm a} \|^2 }
\end{align}
and concludes the proof. 
\end{proof}

\subsection{Proof of Corollary~\ref{lem:normalize}} \label{app:normalize}
\begin{proof}
Suppose that ${\bm x}$ is an NE to the game ${\sf Game}( {\bm a}, {\bm b}, {\bf G} , f ; {\bm p} )$. We note that it must satisfy
\beq
x_i = \max \left\{ 0, \frac{ a_i - p_i }{ 2 b_i } + \frac{1}{2 b_i} \sum_{j=1}^n G_{ij} f( x_i ) \right\} 
\eeq 
From the definition of $\hat{f}(x)$, the above is equivalent to
\beq \notag
\begin{aligned}
& x_i = \max \left\{ 0, \frac{ a_i - p_i }{ 2 b_i } + \frac{1}{2 b_i} \sum_{j=1}^n G_{ij} ( \alpha \hat{f}( x_i ) + f(0) ) \right\} \\
& = \max \left\{ 0, \frac{ a_i + \sum_{j=1}^n G_{ij} f(0) - p_i }{ 2 b_i } + \frac{1}{2 b_i} \sum_{j=1}^n \hat{G}_{ij}  \hat{f}( x_i ) ) \right\}
\end{aligned}
\eeq 
The proof is concluded by observing that the latter is the equilibrium condition for ${\sf Game}( \hat{\bm a}, \hat{\bm b}, \hat{\bf G}, \hat{f} ; {\bm p} )$.
\end{proof}

\subsection{Proof of Lemma~\ref{lem:positive_x}} \label{app:positive_x}

\begin{proof}
Our proof follows the same strategy from \cite[Lemma 7]{candogan2012optimal} with two modifications: (1) the interaction function is nonlinear, (2) the optimal prices may be negative. 

Let ${\bm p}^\star$ be an optimal price that solves \eqref{eq:opt-price} and ${\bm x}^\star$ be the corresponding NE action vector. For the sake of contradiction, we assume that for some $i \in {\cal V}$, $x_i^\star = 0$. We notice that $p_i^\star \geq 0$ since otherwise, from \eqref{eq:ne_eqn}, one must have $x_i^\star > a_i - p_i^\star > 0$.

We construct an alternative price vector ${\bm p}'$ with $p_i' = c \in (0, a_i)$ and for any $j \neq i$,
\[
p_j' = p_j^\star + G_{ji} f \left( \frac{ a_i - c + \sum_{k=1}^n G_{ik} f(x_k^\star) }{ 2 b_i } \right) > p_j^\star
\]
We claim that the action vector ${\bm x}'$ with $x_i' = \frac{ a_i - c + \sum_{k=1}^n G_{ik} f(x_k^\star) }{ 2 b_i } > 0$, $x_j' = x_j^\star$ for any $j \neq i$ is a NE to the game induced by the price vector ${\bm p}'$. From the construction, it is clear that $x_i'$ satisfies \eqref{eq:ne_eqn}. For any $j \neq i$, 
\begin{equation}
    \begin{aligned}
        & T_j( {\bm x}' ) = \max \left\{ 0, \frac{ a_j - p_j' + \sum_{k=1}^n G_{jk} f( x_k' ) }{ 2 b_j } \right\} \\
        & = \max \left\{ 0, \frac{ a_j - p_j' + G_{ji} f( x_i') + \sum_{k \neq i} G_{jk} f( x_k^\star ) }{ 2 b_j } \right\}
    \end{aligned}
\end{equation}
Notice that as $G_{ji} f( x_i') - p_j' = - p_j^\star$, we have 
\begin{align}
    T_j( {\bm x}' ) = \max \left\{ 0, \frac{ a_j - p_j^\star + \sum_{k \neq i} G_{jk} f( x_k^\star ) }{ 2 b_j } \right\} .
\end{align}
The right hand side equals $x_j^\star = x_j'$. As such, ${\bm x}'$ is an NE action vector under ${\bm p}'$. Together, this yields a contradiction to the optimality of ${\bm p}^\star$ since $( {\bm p}' )^\top {\bm x}' > ({\bm p}^\star)^\top {\bm x}^\star$.
\end{proof}

\subsection{Proof of Corollary~\ref{cor:symmetric}} \label{app:symmetric}
\begin{proof}
First observe that under Assumptions~\ref{assu:f_strong}, \ref{assu:B}, there exists $\bar{x} \geq 0$ which satisfies the following condition:
\beq 
\bar{x} = \max \left\{ 0, \frac{ \bar{a} - \bar{p} + \bar{g} f(\bar{x}) }{ 2\bar{b} } \right\}.
\eeq 
The above condition coincides with the equilibrium condition for ${\sf Game}( \bar{a}{\bf 1}, \bar{b}{\bf 1}, {\bf G}, f; \bar{p}{\bf 1} )$ (cf.~\eqref{eq:ne_eqn}). This shows that ${\bm x} = \bar{x} {\bf 1}$ is an NE of the game. The proof is then concluded using the uniqueness of NE under Assumption~\ref{assu:B}. 
\end{proof}
    
\subsection{Proof of Theorem~\ref{the:the3}} \label{app:thm3}
\begin{proof} With a slight abuse of notation, we set $J(x) = J( x {\bf 1} )$ such that $J(\cdot)$ is a one-dimensional function.
Using the Taylor's theorem, we note that there exists $\xi \in ( \bar{x}^\star, \bar{x}_0)$ such that 
\begin{equation} \label{eq:jxbound}
    \begin{aligned}
        J( \bar{x}_0) = J( \bar{x}^\star ) + \frac{1}{2} J''(\xi)( \bar{x}_0 - \bar{x}^\star )^2,
    \end{aligned}
\end{equation}
since $J'(\bar{x}^\star) = 0$.
Notice that for any $x > 0$,
\beq 
\begin{aligned}
\frac{ J''(x) }{ n } = -4 \bar{b} + \bar{g} ( 2 f'(x) + x f''(x) ) \leq -4 \bar{b} + 2 \bar{g}.
\end{aligned}
\eeq 

Since $J'( \bar{x}^\star ) = 0$, we have 
\beq
\begin{aligned}
& \bar{a} - 4 \bar{b} \bar{x}^\star + \bar{g}\left(f( \bar{x}^\star ) + \bar{x}^\star  f'( \bar{x}^\star )\right) = 0.\\
% & \frac{ \bar{a} }{4 \bar{b} } \leq \bar{x}^\star = \frac{ \bar{a} + \bar{g} \left(f'(\bar{x}^\star ) \bar{x}^\star + f( \bar{x}^\star ) \right)}{4 \bar{b} } \leq \frac{1}{1 - \bar{g} / (2 \bar{b}) } \frac{ \bar{a} }{4 \bar{b} }
\end{aligned}
\eeq
It also holds that 
\begin{align}
\bar{a} - 4 \bar{b} \bar{x}_0 + 2 \bar{g} f( \bar{x}_0 ) = 0.
\end{align}
Let \(u(x) = \bar{a} - 4\bar{b}x + \bar{g}\left(xf'(x) + f(x)\right)\) and \(v(x) = \bar{a} - 4\bar{b}x + 2\bar{g}f(x)\). We know \(u(\bar{x}^\star)=0\), \(v(\bar{x}_0) = 0\). Furthermore,
\beq
\begin{aligned}
u'(x) &= -4\bar{b}+\bar{g}(2f'(x) + xf''(x))\\
&\leq -4\bar{b} + 2\bar{g} < 0 \\
v'(x) & = -4\bar{b} + 2\bar{g} f'(x) \leq -4\bar{b} + 2\bar{g} < 0
\end{aligned}
\eeq
Therefore, \(u(x)\), $v(x)$ are decreasing functions with the unique roots given by $\bar{x}^\star, \bar{x}_0$, respectively. 
% and \(\bar{x}^\star\) is the unique solution
% \[
% v'(x) = -4\bar{b} + 2\bar{g} f'(x) \leq -4\bar{b} + 2\bar{g} \leq 0
% \]
% Therefore, \(v(x)\) is decreasing and \(\bar{x}_0\) is the unique solution.

We first prove that $\bar{x}_0 \geq \bar{x}^\star$ by observing
\begin{align}
v(x) - u(x) = \bar{g} ( f(x) - x f'(x) ) \geq 0
\end{align}
where the last inequality is due to the concavity of $f$. The above implies that 
\begin{align}
v( \bar{x}^\star ) \geq u ( \bar{x}^\star ) = 0 = v( \bar{x}_0 )
\end{align}
This implies $\bar{x}^\star \leq \bar{x}_0$ as $v(\cdot)$ is a decreasing function.
% Next, we prove that \(\bar{x}_0 \geq \bar{x}^\star\)
% \[
% \begin{aligned}
% &v(x) - u(x) = \bar{g}\left(f(x) - xf'(x)\right)=\bar{g}{\bf h}(x) > 0\\
% & \ \forall x \in \mathbb{R}_+ \ \  v(x) \geq u(x)\\
% & \ v(\bar{x}^\star) \geq u(\bar{x}^\star) = 0 = v(\bar{x}_0)\\
% & \ \bar{x}_0 \geq \bar{x}^\star \text{ since \(v\) is decreasing.}
% \end{aligned}
% \]
% And since $\bar{b}\bar{x}_0 =  \bar{a} + 2\bar{g}f(\bar{x}_0)
%       \leq  \bar{a} + \bar{g} \bar{x}_0$, we have
%       $\bar{x}_0 < \frac{\bar{a}}{4(1 - \gamma) \bar{b}}$
% Therefore $\frac{ \bar{a} }{ 4 \bar{b} } \leq \bar{x}^\star \leq \bar{x}_0 \leq \frac{1}{1 - \bar{g}/2\bar{b}} \frac{ \bar{a} }{ 4 \bar{b} }$
By noting that $\bar{a} - 4\bar{b}\bar{x}^\star+\bar{g}\left(\bar{x}^\star f'(\bar{x}^\star)+f(\bar{x}^\star)\right) = \bar{a} - 4\bar{b}\bar{x}_0+2\bar{g}f(\bar{x}_0)$, we have
\beq
\begin{aligned}
4\bar{b}(\bar{x}_0- \bar{x}^\star) = & \ \bar{g}\left(2f( \bar{x}_0)-\bar{x}^\star f'(\bar{x}^\star)-f(\bar{x}^\star)\right)\\
 \geq & \ \bar{g}\left( f( \bar{x}_0 )-\bar{x}_0f'(\bar{x}_0) \right) = \bar{g} h(\bar{x}_0),
\end{aligned}
\eeq
where the inequality is due to that $f(x) + f'(x) x$ is increasing in $x$ and $\bar{x}^\star \leq \bar{x}_0$.
As $\frac{ \bar{a} }{ 4 \bar{b} } \leq \bar{x}_0 \leq \frac{ \bar{a} }{ 4 \bar{b} - 2 \bar{g} }$, 
\begin{align}
0 < J( \bar{x}_0 ) = \frac{n}{2} \bar{a} \bar{x}_0 \leq \frac{ n \bar{a}^2 }{ 8 \bar{b} - 4 \bar{g} }.
\end{align}

Substituting the above estimates back into \eqref{eq:jxbound}, we get
\beq 
\begin{aligned}
    \frac{ J( \bar{x}^\star )}{ J( \bar{x}_0 )} - 1 & = - \frac{1}{2} \frac{ J''(\zeta) ( \bar{x}_0 - \bar{x}^\star)^2 }{ J( \bar{x}_0 ) } \\
    & \geq (2 \bar{b} - \bar{g} ) \frac{ n ( \bar{x}_0 - \bar{x}^\star)^2 }{ J( \bar{x}_0 ) } \\
    & \geq \left(  \frac{ \bar{g} }{ \bar{a} }  \left(1 - \frac{ \bar{g} }{ 2 \bar{b} } \right)  h\left(\frac{ \bar{a} }{ 4 \bar{b} } \right) \right)^2.
\end{aligned}
\eeq 
This concludes the proof.
\end{proof}
\bibliographystyle{IEEEtran}
\bibliography{Ref}

\end{document}

%% file: figs_tex/star_pa.tex
\begin{tikzpicture}
\definecolor{darkgray176}{RGB}{176,176,176}
\definecolor{green01270}{RGB}{0,127,0}
\definecolor{lightgray204}{RGB}{204,204,204}

\begin{groupplot}[group style={group name=myplot,group size=2 by 1}]
\nextgroupplot[
legend cell align={left},
legend style={
  fill opacity=0.8,
  draw opacity=1,
  text opacity=1,
  at={(0.5,0.09)},
  anchor=south,
  draw=lightgray204
},
tick align=outside,
tick pos=left,
x grid style={darkgray176},
xlabel={\(\displaystyle \alpha\)},
xmajorgrids, xminorgrids,
xmin=-0.05, xmax=1.05,
xtick style={color=black},
y grid style={darkgray176},
ylabel={\(\displaystyle (PoI)^{-1} = \frac{J({\bm x}_0)}{J({\bm x}^\star)}\)},
ymajorgrids, yminorgrids,
ymin=0.4, ymax=1.03,
ytick style={color=black},
width=6.5cm, 
height=7.5cm,
minor tick num=2 
]
    % \draw[step=0.1cm,gray,very thin] (0.45,0.55) grid (0.9,1);

\addplot [very thick, red]
table {%
0 0.437499999995474
0.01 0.447467334826406
0.02 0.457683858143541
0.03 0.468152269846881
0.04 0.478874781043322
0.05 0.489853044046657
0.06 0.501088076946541
0.07 0.512580182118835
0.08 0.524328859015266
0.09 0.536332710122038
0.1 0.548589341632256
0.11 0.561095257890747
0.12 0.573845750190615
0.13 0.586834780814619
0.14 0.600054862079999
0.15 0.613496932462619
0.16 0.62715022931765
0.17 0.641002161020955
0.18 0.655038179255293
0.19 0.669241653517215
0.2 0.683593749920515
0.21 0.698073317578985
0.22 0.712656784402973
0.23 0.727318066461442
0.24 0.742028493781979
0.25 0.756756756668034
0.26 0.771468876660985
0.27 0.786128206202212
0.28 0.800695461098614
0.29 0.815128790232947
0.3 0.829383886163834
0.31 0.843414140356792
0.32 0.857170846280788
0.33 0.870603452446166
0.34 0.883659866590859
0.35 0.896286811690728
0.36 0.908430232418189
0.37 0.920035749852552
0.38 0.931049159276575
0.39 0.941416966907503
0.4 0.951086956390161
0.41 0.960008777110759
0.42 0.968134542914123
0.43 0.975419430241634
0.44 0.981822261979499
0.45 0.987306064755304
0.46 0.99183858524836
0.47 0.995392753420738
0.48 0.997947080165539
0.49 0.999485978534488
0.5 0.99999999987931
0.51 0.999485978534488
0.52 0.997947080189677
0.53 0.995392753420738
0.54 0.99183858524836
0.55 0.987306064755304
0.56 0.981822262003636
0.57 0.975419430241634
0.58 0.968134542914123
0.59 0.960008777134897
0.6 0.951086956390161
0.61 0.941416966931641
0.62 0.931049159276575
0.63 0.920035749876689
0.64 0.908430232466465
0.65 0.896286811690728
0.66 0.883659866614997
0.67 0.870603452470304
0.68 0.857170846280788
0.69 0.84341414038093
0.7 0.829383886187972
0.71 0.815128790257085
0.72 0.800695461098614
0.73 0.78612820622635
0.74 0.771468876685123
0.75 0.756756756692172
0.76 0.742028493830255
0.77 0.727318066509718
0.78 0.712656784402973
0.79 0.698073317603123
0.8 0.683593749944652
0.81 0.669241653541353
0.82 0.655038179303569
0.83 0.641002161020955
0.84 0.62715022931765
0.85 0.613496932486757
0.86 0.600054862104137
0.87 0.586834780814619
0.88 0.573845750190615
0.89 0.561095257914885
0.9 0.548589341656394
0.91 0.536332710122038
0.92 0.524328859015266
0.93 0.512580182118835
0.94 0.501088076946541
0.95 0.489853044070795
0.96 0.47887478106746
0.97 0.468152269846881
0.98 0.457683858143541
0.99 0.447467334826406
1 0.437499999995474
};
% \addlegendentry{linear interaction}
\addplot [very thick, blue]
table {%
0 0.741946484190034
0.01 0.752492851837965
0.02 0.762951312056385
0.03 0.773307055343185
0.04 0.783545393869663
0.05 0.793651832545729
0.06 0.803612137642251
0.07 0.813412400737801
0.08 0.823039098682031
0.09 0.832479148625409
0.1 0.841719957148925
0.11 0.850749463915184
0.12 0.859556179357381
0.13 0.868129215704732
0.14 0.876458312814509
0.15 0.884533857461802
0.16 0.89234689698666
0.17 0.899889147907414
0.18 0.907152998929365
0.19 0.914131509821795
0.2 0.920818405913916
0.21 0.927208068412146
0.22 0.933295522159062
0.23 0.939076419788487
0.24 0.944547023922055
0.25 0.949704187300369
0.26 0.954545331027146
0.27 0.95906842209704
0.28 0.96327194932827
0.29 0.967154899464072
0.3 0.970716732567345
0.31 0.973957357387111
0.32 0.976877106811674
0.33 0.979476714022092
0.34 0.981757288665372
0.35 0.983720294047126
0.36 0.98536752482488
0.37 0.986701085600853
0.38 0.987723370289318
0.39 0.988437042532924
0.4 0.988845017007836
0.41 0.988950441723801
0.42 0.988756681123212
0.43 0.988267300440883
0.44 0.987486050730993
0.45 0.986416854942399
0.46 0.98506379500145
0.47 0.983431099458945
0.48 0.981523132208831
0.49 0.979344382148851
0.5 0.976899453036276
0.51 0.974193054744779
0.52 0.971229994857628
0.53 0.968015170595818
0.54 0.964553562047727
0.55 0.960850225368285
0.56 0.956910286886338
0.57 0.95273893716266
0.58 0.948341426222207
0.59 0.94372305831027
0.6 0.938889188108529
0.61 0.933845216017122
0.62 0.928596584844187
0.63 0.923148776092465
0.64 0.917507306737488
0.65 0.911677726154896
0.66 0.905665613200643
0.67 0.89947657328926
0.68 0.893116236334463
0.69 0.886590253574229
0.7 0.879904295430657
0.71 0.873064049126831
0.72 0.866075216338017
0.73 0.858943510637709
0.74 0.851674655377351
0.75 0.844274381429339
0.76 0.83674842419559
0.77 0.829102521813077
0.78 0.821342412090169
0.79 0.813473830326339
0.8 0.805502506142028
0.81 0.797434161282157
0.82 0.789274505905767
0.83 0.781029236250786
0.84 0.772704030923405
0.85 0.764304547977038
0.86 0.755836421485074
0.87 0.74730525758602
0.88 0.73871663117433
0.89 0.730076081817392
0.9 0.721389110182933
0.91 0.712661173428659
0.92 0.703897681390142
0.93 0.695103992051499
0.94 0.686285407437458
0.95 0.677447168752569
0.96 0.668594452012859
0.97 0.659732363514788
0.98 0.650865934852104
0.99 0.642000118613619
1 0.633139783416954
};
% \addlegendentry{$f_1(x)$}
\addplot [very thick, green01270]
table {%
0 0.959296889010653
0.01 0.961689688145903
0.02 0.963983392324506
0.03 0.966180032475237
0.04 0.968281577188483
0.05 0.970289933635892
0.06 0.972206950482844
0.07 0.974034419239824
0.08 0.975774075686567
0.09 0.97742760252443
0.1 0.978996630611427
0.11 0.980482740812103
0.12 0.981887465839546
0.13 0.983212291866721
0.14 0.984458659988013
0.15 0.985627968258415
0.16 0.986721572792415
0.17 0.987740789540028
0.18 0.988686895635844
0.19 0.98956113092708
0.2 0.99036469906734
0.21 0.991098769291908
0.22 0.991764477200322
0.23 0.992362926464936
0.24 0.992895189363959
0.25 0.993362308550455
0.26 0.993765298075446
0.27 0.994105144166184
0.28 0.994382806247126
0.29 0.994599218270471
0.3 0.994755289303828
0.31 0.994851904733697
0.32 0.994889926789555
0.33 0.994870195867521
0.34 0.99479353074525
0.35 0.994660729902417
0.36 0.994472571979372
0.37 0.994229816481873
0.38 0.993933204788004
0.39 0.993583460237913
0.4 0.993181289507939
0.41 0.992727382573364
0.42 0.992222413591779
0.43 0.991667041598221
0.44 0.991061910773739
0.45 0.990407651139772
0.46 0.989704879007132
0.47 0.988954197607587
0.48 0.988156197358244
0.49 0.987311456429745
0.5 0.986420541070518
0.51 0.985484006354005
0.52 0.984502396137526
0.53 0.983476243626103
0.54 0.982406072056695
0.55 0.98129239453367
0.56 0.98013571464792
0.57 0.978936526979898
0.58 0.977695317052686
0.59 0.976412562007246
0.6 0.975088730562045
0.61 0.973724283565317
0.62 0.972319674189178
0.63 0.970875348059979
0.64 0.969391743936695
0.65 0.967869293475502
0.66 0.966308421665717
0.67 0.964709547136345
0.68 0.963073082346811
0.69 0.96139943365526
0.7 0.95968900186485
0.71 0.957942182054758
0.72 0.956159363942614
0.73 0.954340932192487
0.74 0.952487266360421
0.75 0.950598741315871
0.76 0.948675727188865
0.77 0.946718589732945
0.78 0.944727690269123
0.79 0.942703386100082
0.8 0.940646030406537
0.81 0.938555972656803
0.82 0.936433558370106
0.83 0.934279129830788
0.84 0.932093025621101
0.85 0.929875581034266
0.86 0.927627128188747
0.87 0.925347996141215
0.88 0.923038511202533
0.89 0.92069899643784
0.9 0.91832977245036
0.91 0.915931157133522
0.92 0.913503465849876
0.93 0.911047011594182
0.94 0.90856210489797
0.95 0.906049054264745
0.96 0.903508165775396
0.97 0.900939743770498
0.98 0.898344090451059
0.99 0.895721506043893
1 0.893072289133711
};
% \addlegendentry{$f_2(x)$}

\nextgroupplot[
legend cell align={left},
legend style={
  fill opacity=0.8,
  draw opacity=1,
  text opacity=1,
  at={(0.03,0.03)},
  anchor=south west,
  draw=lightgray204,
},
tick align=outside,
tick pos=left,
x grid style={darkgray176},
xlabel={\(\displaystyle \alpha\)},
xmajorgrids, xminorgrids,
xmin=-0.05, xmax=1.05,
xtick style={color=black},
y grid style={darkgray176},
ymajorgrids, yminorgrids,
ymin=0.7, ymax=1.03,
ytick style={color=black},
width=6.5cm, 
height=7.5cm,
minor tick num = 2
]
\addplot [very thick, red]
table {%
0 0.821462598009049
0.01 0.826562768656991
0.02 0.831669654318631
0.03 0.836780539671634
0.04 0.841892555018127
0.05 0.847002672761526
0.06 0.852107704336337
0.07 0.857204297520478
0.08 0.862288934200572
0.09 0.86735792871103
0.1 0.872407426736824
0.11 0.877433404910846
0.12 0.882431671146088
0.13 0.887397865662398
0.14 0.892327463211116
0.15 0.897215775732552
0.16 0.902057956422744
0.17 0.906849004766561
0.18 0.911583772597574
0.19 0.916256971607448
0.2 0.92086318165058
0.21 0.925396860780113
0.22 0.929852356330848
0.23 0.934223917532215
0.24 0.93850570916813
0.25 0.942691826927984
0.26 0.94677631385484
0.27 0.950753178333768
0.28 0.954616412935811
0.29 0.958360014993331
0.3 0.961978007678675
0.31 0.965464462763895
0.32 0.96881352365226
0.33 0.972019429466716
0.34 0.975076539923511
0.35 0.977979360074968
0.36 0.980722565958206
0.37 0.983301029670091
0.38 0.985709844744181
0.39 0.987944350943816
0.4 0.990000158259989
0.41 0.99187317045623
0.42 0.993559607113628
0.43 0.99505602473969
0.44 0.996359336045153
0.45 0.997466827710859
0.46 0.998376176010528
0.47 0.999085460581344
0.48 0.999593175728321
0.49 0.999898239463769
0.5 0.999999999939603
0.51 0.999898239463769
0.52 0.999593175728321
0.53 0.999085460581344
0.54 0.998376176010528
0.55 0.997466827710859
0.56 0.996359336065286
0.57 0.995056024759823
0.58 0.993559607113628
0.59 0.99187317045623
0.6 0.990000158290188
0.61 0.987944350943816
0.62 0.98570984477438
0.63 0.983301029690224
0.64 0.980722565958206
0.65 0.977979360105167
0.66 0.975076539923511
0.67 0.972019429496914
0.68 0.968813523662327
0.69 0.965464462794094
0.7 0.961978007708873
0.71 0.958360015013463
0.72 0.954616412976076
0.73 0.950753178363967
0.74 0.946776313885038
0.75 0.94269182693805
0.76 0.938505709188263
0.77 0.934223917552347
0.78 0.929852356381179
0.79 0.925396860810312
0.8 0.920863181680778
0.81 0.916256971637647
0.82 0.911583772627773
0.83 0.906849004776628
0.84 0.902057956442876
0.85 0.897215775762751
0.86 0.892327463261447
0.87 0.887397865712729
0.88 0.882431671156154
0.89 0.877433404930979
0.9 0.872407426756957
0.91 0.867357928741228
0.92 0.862288934230771
0.93 0.857204297550676
0.94 0.852107704366536
0.95 0.847002672791724
0.96 0.841892555048325
0.97 0.836780539691766
0.98 0.831669654328698
0.99 0.826562768656991
1 0.821462598009049
};
\addlegendentry{linear interaction}
\addplot [very thick, blue]
table {%
0 0.900741397573379
0.01 0.904808141488169
0.02 0.908819908239412
0.03 0.912773904994703
0.04 0.91666736115125
0.05 0.920497534588696
0.06 0.924261717856081
0.07 0.92795724417804
0.08 0.931581493591371
0.09 0.935131898441702
0.1 0.93860594916527
0.11 0.942001199821167
0.12 0.94531527291889
0.13 0.948545864601583
0.14 0.95169074894228
0.15 0.954747782506189
0.16 0.957714908096632
0.17 0.960590158199708
0.18 0.963371658579268
0.19 0.96605763068166
0.2 0.968646394445399
0.21 0.971136370299396
0.22 0.973526080701849
0.23 0.975814151749595
0.24 0.977999314080045
0.25 0.980080403336959
0.26 0.982056360731582
0.27 0.983926232855128
0.28 0.985689171323289
0.29 0.987344432250591
0.3 0.988891375113535
0.31 0.990329461835221
0.32 0.991658255169881
0.33 0.992877417209521
0.34 0.993986707448577
0.35 0.994985980737769
0.36 0.995875185263116
0.37 0.996654360126978
0.38 0.997323632887518
0.39 0.997883217101324
0.4 0.998333409602752
0.41 0.998674587931135
0.42 0.99890720739566
0.43 0.999031798454594
0.44 0.999048963799162
0.45 0.998959375541402
0.46 0.99876377236185
0.47 0.998462956650539
0.48 0.99805779177403
0.49 0.997549199169834
0.5 0.996938155644309
0.51 0.996225690619721
0.52 0.995412883474124
0.53 0.994500860896616
0.54 0.993490794291053
0.55 0.992383897239628
0.56 0.991181423091514
0.57 0.98988466244587
0.58 0.988494941046641
0.59 0.987013617163639
0.6 0.985442079638541
0.61 0.983781745615233
0.62 0.982034058527766
0.63 0.980200486046974
0.64 0.9782825180986
0.65 0.976281665003367
0.66 0.974199455743587
0.67 0.972037435991021
0.68 0.969797166587576
0.69 0.967480221777151
0.7 0.965088187781717
0.71 0.962622661060888
0.72 0.96008524699115
0.73 0.957477558375066
0.74 0.954801214081839
0.75 0.95205783776935
0.76 0.949249056551538
0.77 0.946376499727947
0.78 0.943441797783149
0.79 0.940446581054411
0.8 0.937392478696819
0.81 0.934281117646655
0.82 0.931114121589908
0.83 0.927893110022923
0.84 0.924619697147229
0.85 0.921295491213817
0.86 0.917922093404703
0.87 0.914501097074009
0.88 0.911034087031459
0.89 0.907522638552639
0.9 0.903968316838836
0.91 0.900372676103606
0.92 0.896737258989538
0.93 0.893063595870271
0.94 0.889353204196706
0.95 0.885607587804798
0.96 0.88182823644185
0.97 0.878016625097887
0.98 0.874174213490766
0.99 0.870302445514179
1 0.866402748791515
};
\addlegendentry{$f_1(x)$}
\addplot [very thick, green01270]
table {%
0 0.98138241885056
0.01 0.982528664413657
0.02 0.983627539022482
0.03 0.984679926471639
0.04 0.98568669706055
0.05 0.986648706006555
0.06 0.987566792982562
0.07 0.988441781361063
0.08 0.989274477246516
0.09 0.990065669629499
0.1 0.990816129781038
0.11 0.991526611301705
0.12 0.992197849716754
0.13 0.992830562886808
0.14 0.993425450910669
0.15 0.993983196255819
0.16 0.994504464035586
0.17 0.994989902166602
0.18 0.995440141765731
0.19 0.995855797214523
0.2 0.996237466809407
0.21 0.996585732865932
0.22 0.996901162195331
0.23 0.99718430639531
0.24 0.997435702339632
0.25 0.99765587242846
0.26 0.997845325051281
0.27 0.998004554982933
0.28 0.998134043726489
0.29 0.998234259855919
0.3 0.998305659517738
0.31 0.998348686693622
0.32 0.998363773568821
0.33 0.998351340953575
0.34 0.998311798558218
0.35 0.99824554541401
0.36 0.998152970174319
0.37 0.998034451349435
0.38 0.997890357806085
0.39 0.997721048908894
0.4 0.997526874926793
0.41 0.99730817729323
0.42 0.997065288839946
0.43 0.996798534149408
0.44 0.996508229788185
0.45 0.996194684566164
0.46 0.995858199743287
0.47 0.995499069380807
0.48 0.995117580468131
0.49 0.994714013168257
0.5 0.994288641102655
0.51 0.99384173150357
0.52 0.993373545419625
0.53 0.992884337960257
0.54 0.992374358408275
0.55 0.991843850503627
0.56 0.991293052542344
0.57 0.990722197607479
0.58 0.990131513667484
0.59 0.989521223794038
0.6 0.988891546338311
0.61 0.988242694991403
0.62 0.987574879025831
0.63 0.986888303394207
0.64 0.986183168840906
0.65 0.985459672090753
0.66 0.984718005921433
0.67 0.983958359326049
0.68 0.983180917624268
0.69 0.982385862559017
0.7 0.981573372433964
0.71 0.980743622196592
0.72 0.979896783561601
0.73 0.979033025082607
0.74 0.978152512286143
0.75 0.97725540775758
0.76 0.976341871222649
0.77 0.975412059605745
0.78 0.97446612720384
0.79 0.973504225653832
0.8 0.97252650409077
0.81 0.971533109273286
0.82 0.970524185507575
0.83 0.969499874822333
0.84 0.968460317038322
0.85 0.967405649786106
0.86 0.966336008690724
0.87 0.965251527248737
0.88 0.964152337038882
0.89 0.96303856773801
0.9 0.961910347180373
0.91 0.960767801423298
0.92 0.959611054752377
0.93 0.958440229764945
0.94 0.957255447512528
0.95 0.956056827353843
0.96 0.954844487226347
0.97 0.953618543510756
0.98 0.952379111173893
0.99 0.951126303789237
1 0.949860233587277
};
\addlegendentry{$f_2(x)$}

\end{groupplot}
\end{tikzpicture}%

%% file: figs_tex/bound_all.tex
\begin{tikzpicture}
\definecolor{darkgray176}{RGB}{176,176,176}
\definecolor{green01270}{RGB}{0,127,0}
\definecolor{lightgray204}{RGB}{204,204,204}

\begin{groupplot}[group style={group name=myplot,group size=2 by 1}]
\nextgroupplot[
legend cell align={left},
legend style={fill opacity=0.8, draw opacity=1, text opacity=1, draw=lightgray204},
tick align=outside,
tick pos=left,
x grid style={darkgray176},
xlabel={\(\displaystyle \bar{g}\)},
xmajorgrids,
xmin=-0.0395, xmax=1.0495,
xtick style={color=black},
y grid style={darkgray176},
ylabel={Value},
ymajorgrids, yminorgrids,
ymin=3e-06, ymax=0.008,
ymode=log,
ytick style={color=black},
width=6.5cm, 
height=7.5cm,
minor tick num = 1
]
\addplot [very thick, red]
table {%
0.01 0.0001316466
0.02 0.0002840427
0.03 0.0004149734
0.04 0.0005226857
0.05 0.0006099769
0.06 0.0006801395
0.07 0.0007361348
0.08 0.000780452
0.09 0.0008151377
0.1 0.0008418639
0.11 0.0008619963
0.12 0.0008766534
0.13 0.0008867543
0.14 0.0008930575
0.15 0.0008961916
0.16 0.0008966796
0.17 0.0008949581
0.18 0.0008913934
0.19 0.0008862929
0.2 0.0008799164
0.21 0.0008724831
0.22 0.000864179
0.23 0.000855162
0.24 0.0008455663
0.25 0.0008355064
0.26 0.0008250798
0.27 0.0008143699
0.28 0.0008034477
0.29 0.000792374
0.3 0.0007812007
0.31 0.0007699722
0.32 0.0007587262
0.33 0.000747495
0.34 0.0007363059
0.35 0.0007251821
0.36 0.0007141432
0.37 0.0007032059
0.38 0.0006923838
0.39 0.0006816885
0.4 0.0006711296
0.41 0.0006607147
0.42 0.0006504501
0.43 0.0006403408
0.44 0.0006303905
0.45 0.000620602
0.46 0.0006109773
0.47 0.0006015177
0.48 0.0005922237
0.49 0.0005830954
0.5 0.0005741324
0.51 0.0005653338
0.52 0.0005566985
0.53 0.0005482249
0.54 0.0005399115
0.55 0.0005317563
0.56 0.0005237572
0.57 0.000515912
0.58 0.0005082183
0.59 0.0005006737
0.6 0.0004932757
0.61 0.0004860216
0.62 0.0004789089
0.63 0.0004719348
0.64 0.0004650969
0.65 0.0004583923
0.66 0.0004518184
0.67 0.0004453725
0.68 0.000439052
0.69 0.0004328543
0.7 0.0004267767
0.71 0.0004208167
0.72 0.0004149718
0.73 0.0004092394
0.74 0.0004036171
0.75 0.0003981025
0.76 0.0003926932
0.77 0.0003873868
0.78 0.0003821811
0.79 0.0003770738
0.8 0.0003720627
0.81 0.0003671456
0.82 0.0003623206
0.83 0.0003575854
0.84 0.000352938
0.85 0.0003483766
0.86 0.0003438991
0.87 0.0003395036
0.88 0.0003351885
0.89 0.0003309517
0.9 0.0003267917
0.91 0.0003227066
0.92 0.0003186948
0.93 0.0003147546
0.94 0.0003108846
0.95 0.000307083
0.96 0.0003033485
0.97 0.0002996795
0.98 0.0002960745
0.99 0.0002925322
1 0.0002890512
};
\addlegendentry{$\frac{J({\bm x}^*)}{J({\bm x}_0)}-1$}
\addplot [very thick, blue]
table {%
0.01 3.754762e-06
0.02 1.693893e-05
0.03 3.384718e-05
0.04 5.053095e-05
0.05 6.539505e-05
0.06 7.797166e-05
0.07 8.827609e-05
0.08 9.651355e-05
0.09 0.0001029483
0.1 0.0001078465
0.11 0.0001114528
0.12 0.0001139823
0.13 0.0001156197
0.14 0.0001165216
0.15 0.0001168196
0.16 0.000116624
0.17 0.0001160265
0.18 0.0001151038
0.19 0.0001139196
0.2 0.0001125269
0.21 0.0001109699
0.22 0.0001092851
0.23 0.0001075031
0.24 0.0001056492
0.25 0.0001037443
0.26 0.0001018059
0.27 9.984845e-05
0.28 9.78837e-05
0.29 9.59215e-05
0.3 9.396983e-05
0.31 9.203525e-05
0.32 9.012304e-05
0.33 8.823744e-05
0.34 8.638184e-05
0.35 8.455886e-05
0.36 8.277054e-05
0.37 8.101839e-05
0.38 7.930348e-05
0.39 7.762652e-05
0.4 7.598794e-05
0.41 7.438789e-05
0.42 7.282635e-05
0.43 7.130308e-05
0.44 6.981774e-05
0.45 6.836986e-05
0.46 6.695889e-05
0.47 6.558417e-05
0.48 6.424504e-05
0.49 6.294073e-05
0.5 6.167049e-05
0.51 6.043352e-05
0.52 5.922901e-05
0.53 5.805613e-05
0.54 5.691407e-05
0.55 5.5802e-05
0.56 5.47191e-05
0.57 5.366455e-05
0.58 5.263757e-05
0.59 5.163736e-05
0.6 5.066314e-05
0.61 4.971416e-05
0.62 4.878968e-05
0.63 4.788897e-05
0.64 4.701133e-05
0.65 4.615606e-05
0.66 4.532251e-05
0.67 4.451001e-05
0.68 4.371793e-05
0.69 4.294567e-05
0.7 4.219262e-05
0.71 4.145821e-05
0.72 4.074188e-05
0.73 4.004309e-05
0.74 3.936131e-05
0.75 3.869604e-05
0.76 3.804678e-05
0.77 3.741305e-05
0.78 3.67944e-05
0.79 3.619039e-05
0.8 3.560057e-05
0.81 3.502454e-05
0.82 3.446189e-05
0.83 3.391222e-05
0.84 3.337518e-05
0.85 3.285038e-05
0.86 3.233748e-05
0.87 3.183614e-05
0.88 3.134603e-05
0.89 3.086682e-05
0.9 3.039821e-05
0.91 2.993991e-05
0.92 2.949162e-05
0.93 2.905307e-05
0.94 2.862398e-05
0.95 2.82041e-05
0.96 2.779317e-05
0.97 2.739095e-05
0.98 2.69972e-05
0.99 2.661169e-05
1 2.62342e-05
};
\addlegendentry{Theorem~\ref{the:the2}}

\nextgroupplot[
legend cell align={left},
legend style={
  fill opacity=0.8,
  draw opacity=1,
  text opacity=1,
  at={(0.03,0.97)},
  anchor=north west,
  draw=lightgray204
},
tick align=outside,
tick pos=left,
x grid style={darkgray176},
xlabel={\(\displaystyle \bar{g}\)},
xmajorgrids,
xmin=-0.0133337, xmax=1.0133337,
xtick style={color=black},
y grid style={darkgray176},
ymajorgrids, yminorgrids,
ymode=log,
ymin=1e-6, ymax=0.01,
ytick style={color=black},
width=6.5cm, 
height=7.5cm,
minor tick num = 1
]
\addplot [very thick, red]
table {%
0.033333 1.500411e-06
0.066667 6.231082e-06
0.1 1.455739e-05
0.133333 2.687411e-05
0.166667 4.36074e-05
0.2 6.521672e-05
0.233333 9.219721e-05
0.266667 0.0001250807
0.3 0.0001644393
0.333333 0.0002108864
0.366667 0.0002650797
0.4 0.0003277227
0.433333 0.000399568
0.466667 0.0004814193
0.5 0.0005741326
0.533333 0.0006786199
0.566667 0.0007958524
0.6 0.0009268593
0.633333 0.001072734
0.666667 0.001234634
0.7 0.001413782
0.733333 0.001611473
0.766667 0.001829068
0.8 0.002068002
0.833333 0.002329784
0.866667 0.002615998
0.9 0.002928301
0.933333 0.003268429
0.966667 0.003638194
};
\addlegendentry{$\frac{J(\bar{x}^\star {\bf 1} )}{J(\bar{x}_0 {\bf 1} )}-1$}
\addplot [very thick, blue]
table {%
0.033333 1.397501e-06
0.066667 5.40212e-06
0.1 1.173925e-05
0.133333 2.014393e-05
0.166667 3.036083e-05
0.2 4.214424e-05
0.233333 5.525811e-05
0.266667 6.947601e-05
0.3 8.458115e-05
0.333333 0.0001003664
0.366667 0.0001166342
0.4 0.0001331966
0.433333 0.0001498755
0.466667 0.0001665022
0.5 0.0001829177
0.533333 0.0001989727
0.566667 0.0002145275
0.6 0.000229452
0.633333 0.0002436257
0.666667 0.0002569379
0.7 0.0002692874
0.733333 0.0002805826
0.766667 0.0002907417
0.8 0.0002996924
0.833333 0.000307372
0.866667 0.0003137276
0.9 0.0003187158
0.933333 0.0003223029
0.966667 0.0003244648
};
\addlegendentry{Theorem~\ref{the:the3}}
\end{groupplot}
\end{tikzpicture}%

%% file: main.bbl
\begin{thebibliography}{10}
\providecommand{\url}[1]{#1}
\csname url@samestyle\endcsname
\providecommand{\newblock}{\relax}
\providecommand{\bibinfo}[2]{#2}
\providecommand{\BIBentrySTDinterwordspacing}{\spaceskip=0pt\relax}
\providecommand{\BIBentryALTinterwordstretchfactor}{4}
\providecommand{\BIBentryALTinterwordspacing}{\spaceskip=\fontdimen2\font plus
\BIBentryALTinterwordstretchfactor\fontdimen3\font minus \fontdimen4\font\relax}
\providecommand{\BIBforeignlanguage}[2]{{%
\expandafter\ifx\csname l@#1\endcsname\relax
\typeout{** WARNING: IEEEtran.bst: No hyphenation pattern has been}%
\typeout{** loaded for the language `#1'. Using the pattern for}%
\typeout{** the default language instead.}%
\else
\language=\csname l@#1\endcsname
\fi
#2}}
\providecommand{\BIBdecl}{\relax}
\BIBdecl

\bibitem{jackson2015games}
M.~O. Jackson and Y.~Zenou, ``Games on networks,'' in \emph{Handbook of game theory with economic applications}, 2015, vol.~4, pp. 95--163.

\bibitem{zhou2016game}
Z.~Zhou, B.~Yolken, R.~A. Miura-Ko, and N.~Bambos, ``A game-theoretical formulation of influence networks,'' in \emph{ACC}, 2016.

\bibitem{parise2017sensitivity}
F.~Parise and A.~Ozdaglar, ``Sensitivity analysis for network aggregative games,'' in \emph{CDC}.\hskip 1em plus 0.5em minus 0.4em\relax IEEE, 2017, pp. 3200--3205.

\bibitem{ballester2006s}
C.~Ballester, A.~Calv{\'o}-Armengol, and Y.~Zenou, ``Who's who in networks. wanted: The key player,'' \emph{Econometrica}, vol.~74, no.~5, pp. 1403--1417, 2006.

\bibitem{candogan2012optimal}
O.~Candogan, K.~Bimpikis, and A.~Ozdaglar, ``Optimal pricing in networks with externalitie,'' \emph{Operation Research}, vol.~60, no.~4, pp. 883--905, 2012.

\bibitem{bramoulle2014strategic}
Y.~Bramoull{\'e}, R.~Kranton, and M.~D'amours, ``Strategic interaction and networks,'' \emph{American Economic Review}, vol. 104, no.~3, 2014.

\bibitem{demange2017optimal}
G.~Demange, ``Optimal targeting strategies in a network under complementarities,'' \emph{Games and Economic Behavior}, vol. 105, 2017.

\bibitem{parise2021analysis}
F.~Parise and A.~Ozdaglar, ``Analysis and interventions in large network games,'' \emph{Annual Review of Control, Robotics, and Autonomous Systems}, vol.~4, pp. 455--486, 2021.

\bibitem{galeotti2009influencing}
A.~Galeotti and S.~Goyal, ``Influencing the influencers: a theory of strategic diffusion,'' \emph{The RAND Journal of Economics}, vol.~40, no.~3, pp. 509--532, 2009.

\bibitem{galeotti2020targeting}
A.~Galeotti, B.~Golub, and S.~Goyal, ``Targeting interventions in networks,'' \emph{Econometrica}, vol.~88, no.~6, pp. 2445--2471, 2020.

\bibitem{maheshwari2022inducing}
C.~Maheshwari, K.~Kulkarni, M.~Wu, and S.~S. Sastry, ``Inducing social optimality in games via adaptive incentive design,'' in \emph{CDC}, 2022.

\bibitem{ata2023latent}
B.~Ata, A.~Belloni, and O.~Candogan, ``Latent agents in networks: Estimation and targeting,'' \emph{Operations Research}, 2023.

\bibitem{belhaj2014network}
M.~Belhaj, Y.~Bramoull{\'e}, and F.~Dero{\"\i}an, ``Network games under strategic complementarities,'' \emph{Games and Economic Behavior}, vol.~88, pp. 310--319, 2014.

\bibitem{cooley2007desegregation}
J.~Cooley, ``Desegregation and the achievement gap: Do diverse peers help?'' \emph{Unpublished manuscript, University of Wisconsin-Madison}, 2007.

\bibitem{naghizadeh2017on}
P.~Naghizadeh and M.~Liu, ``On the uniqueness and stability of equilibria of network games,'' in \emph{Allerton}, 2017.

\bibitem{shende2021network}
A.~Shende, D.~Vasal, and S.~Vishwanath, ``Network design for social welfare,'' in \emph{CISS}, 2021, pp. 1--6.

\bibitem{wang2023network}
X.~Wang, C.-Y. Yau, and H.~T. Wai, ``Network effects in performative prediction games,'' in \emph{ICML}, 2023.

\bibitem{ebrahimi2023united}
R.~Ebrahimi and P.~Naghizadeh, ``United we fall: On the nash equilibria of multiplex network games,'' in \emph{Allerton}.\hskip 1em plus 0.5em minus 0.4em\relax IEEE, 2023, pp. 1--8.

\bibitem{scutari2014real}
G.~Scutari, F.~Facchinei, J.-S. Pang, and D.~P. Palomar, ``Real and complex monotone communication games,'' \emph{IEEE Transactions on Information Theory}, vol.~60, no.~7, pp. 4197--4231, 2014.

\bibitem{parise2019variational}
F.~Parise and A.~Ozdaglar, ``A variational inequality framework for network games: Existence, uniqueness, convergence and sensitivity analysis,'' \emph{Games and Economic Behavior}, vol. 114, pp. 47--82, 2019.

\bibitem{salehisadaghiani2016distributed}
F.~Salehisadaghiani and L.~Pavel, ``Distributed nash equilibrium seeking: A gossip-based algorithm,'' \emph{Automatica}, vol.~72, pp. 209--216, 2016.

\bibitem{liang2017distributed}
S.~Liang, P.~Yi, and Y.~Hong, ``Distributed nash equilibrium seeking for aggregative games with coupled constraints,'' \emph{Automatica}, vol.~85, pp. 179--185, 2017.

\bibitem{tatarenko2020geometric}
T.~Tatarenko, W.~Shi, and A.~Nedi{\'c}, ``Geometric convergence of gradient play algorithms for distributed nash equilibrium seeking,'' \emph{IEEE Transactions on Automatic Control}, vol.~66, no.~11, pp. 5342--5353, 2020.

\bibitem{bacsar1998dynamic}
T.~Ba{\c{s}}ar and G.~J. Olsder, \emph{Dynamic noncooperative game theory}, 1998.

\bibitem{briest2012stackelberg}
P.~Briest, M.~Hoefer, and P.~Krysta, ``Stackelberg network pricing games,'' \emph{Algorithmica}, vol.~62, no.~3, pp. 733--753, 2012.

\bibitem{bohnlein2021revenue}
T.~B{\"o}hnlein, S.~Kratsch, and O.~Schaudt, ``Revenue maximization in stackelberg pricing games: Beyond the combinatorial setting,'' \emph{Mathematical Programming}, vol. 187, pp. 653--695, 2021.

\bibitem{beck2017first}
A.~Beck, \emph{First-order methods in optimization}.\hskip 1em plus 0.5em minus 0.4em\relax SIAM, 2017.

\end{thebibliography}
